\documentclass[10pt]{amsart}
\usepackage{amscd}
\usepackage{amssymb}
\usepackage[all]{xy}
\date{20 July 2008}
\title[Rigid Dualizing Complexes]
{Rigid Dualizing Complexes via Differential Graded Algebras (Survey)}

\author{Amnon Yekutieli}
\address{Department of  Mathematics
Ben Gurion University,
Be'er Sheva 84105,
Israel}
\email{amyekut@math.bgu.ac.il}
\thanks{{\em Mathematics Subject Classification} 2000.
Primary: 18E30; Secondary: 18G10, 16E45, 18G15.}
\keywords{commutative rings, DG algebras,
derived categories, rigid complexes.}
\thanks{This research was supported by the US-Israel Binational
Science Foundation.}

\newtheorem{thm}[equation]{Theorem}
\newtheorem{cor}[equation]{Corollary}

\theoremstyle{definition}
\newtheorem{dfn}[equation]{Definition}
\newtheorem{rem}[equation]{Remark}
\newtheorem{exa}[equation]{Example}

\numberwithin{equation}{section}

\newcommand{\iso}{\stackrel{\simeq}{\rightarrow}}

\newcommand{\xar}{\xrightarrow}
\newcommand{\opn}{\operatorname}
\newcommand{\cat}[1]{\operatorname{\mathsf{#1}}}

\newcommand{\rmitem}[1]{\item[\text{\textup{(#1)}}]}
\newcommand{\mfrak}[1]{\mathfrak{#1}}
\newcommand{\mcal}[1]{\mathcal{#1}}
\newcommand{\msf}[1]{\mathsf{#1}}

\newcommand{\mrm}[1]{\mathrm{#1}}
\newcommand{\mbb}[1]{\mathbb{#1}}

\newcommand{\tup}[1]{\textup{#1}}
\newcommand{\bsym}[1]{\boldsymbol{#1}}

\newcommand{\boplus}{\bigoplus\nolimits}

\DeclareMathSymbol{\mathbbk}{\mathord}{AMSb}{"7C}

\newcommand{\K}{\mbb{K}}

\newcommand{\til}[1]{\tilde{#1}}

\begin{document}

\begin{abstract}
In this article we survey recent results on rigid dualizing 
complexes over commutative algebras. We begin by recalling what are 
dualizing complexes. Next we define rigid complexes, and explain 
their functorial properties. Due to the possible 
presence of torsion, we must use differential graded algebras in the 
constructions. We then discuss rigid dualizing complexes. Finally we 
show how rigid complexes can be used to understand Cohen-Macaulay 
homomorphisms and relative dualizing sheaves.
\end{abstract}

\maketitle

\setcounter{section}{-1}
\section{Introduction}

This short article is based on a lecture I gave at the ``Workshop on 
Triangulated Categories'', Leeds, August 2006. It is a survey of 
recent results on rigid dualizing complexes over commutative rings. 
Most of these results are joint work of mine with James Zhang. The idea of 
rigid dualizing complex is due to Michel Van den Bergh. 

By default all rings considered in this article are {\em commutative}.
We begin by recalling the notion of {\em dualizing complex} 
over a noetherian ring $A$. Next let $B$ be a noetherian $A$-algebra. 
We define what is a {\em rigid complex} 
of $B$-modules relative to $A$. In making this definition 
we must use differential graded algebras (when $B$ is not flat over $A$). 
The functorial properties of rigid complexes are explained. We then 
discuss {\em rigid dualizing complexes}, which by definition are 
complexes that are both rigid and dualizing. Finally we show how rigid complexes 
can be used to understand Cohen-Macaulay homomorphisms and relative dualizing 
sheaves. 

I wish to thank my collaborator James Zhang. Thanks also to Luchezar 
Avramov, Srikanth Iyengar and Joseph Lipman for discussions 
regarding the material in Section 5.

\section{Dualizing Complexes: Overview}

Let $A$ be a noetherian ring. Denote by 
$\msf{D}^{\mrm{b}}_{\mrm{f}}(\cat{Mod} A)$
the derived category of bounded complexes of $A$-modules 
with finitely generated cohomology modules. 

\begin{dfn} (Grothendieck \cite{RD}) \label{dfn1}
A {\em dualizing complex} over $A$ is a complex 
$R \in \msf{D}^{\mrm{b}}_{\mrm{f}}(\cat{Mod} A)$
satisfying the two conditions:
\begin{enumerate}
\rmitem{i} $R$ has finite injective dimension.
\rmitem{ii} The canonical morphism
$A \to \opn{RHom}_{A}(R, R)$ 
is an isomorphism.
\end{enumerate}
\end{dfn}

Condition (i) means that there is an integer $d$ such that
$\opn{Ext}^i_A(M, R) = 0$ for all $i > d$ and all modules $M$.

Recall that a noetherian ring $\K$ is called {\em regular} if all 
its local rings $\K_{\mfrak{p}}$, $\mfrak{p} \in \opn{Spec} \K$, are 
regular local rings. 

\begin{exa} \label{exa0}
If $\K$ is a regular noetherian ring of finite Krull dimension 
(say a field, or the ring of integers $\mbb{Z}$) then 
\[ R := \K \in \msf{D}^{\mrm{b}}_{\mrm{f}}(\cat{Mod} \K) \] 
is a dualizing complex over $\K$. 
\end{exa}

Dualizing complexes over commutative rings are part of 
Grothendieck's duality theory in algebraic geometry, which 
was developed in \cite{RD}. This duality theory deals with 
dualizing complexes on schemes and relations between them. See Remark
\ref{rem1}.

In Section 4 we explain a new approach to dualizing 
complexes over commutative rings, due to James Zhang 
and the author (see \cite{YZ4} and \cite{YZ5}).
Specifically, we discuss existence and 
uniqueness of {\em rigid dualizing complexes}. 

A dualizing complex $R$ has many automorphisms; indeed, its group of automorphisms 
in $\msf{D}(\cat{Mod} A)$ is the group $A^{\times}$ of invertible 
elements. The purpose of rigidity is to eliminate automorphisms, and to
make dualizing complexes functorial. See Theorem \ref{thm4}.

In a sequel paper \cite{Ye2} we use the technique of
{\em perverse coherent sheaves} to construct 
rigid dualizing complexes on schemes, 
and we reproduce almost all of the geometric 
Grothendieck duality theory. 

Related work in noncommutative algebraic geometry (where 
rigid dualizing complexes were first introduced) can be found in 
\cite{VdB, YZ1, YZ2, YZ3}.

\section{Rigid Complexes and DG Algebras}

Let me start with a discussion of rigidity for algebras over a 
field. Suppose $\K$ is a field, $B$ is a $\K$-algebra, and
$M \in \msf{D}(\cat{Mod} B)$. 

According to Van den Bergh 
\cite{VdB} a {\em rigidifying isomorphism} for $M$ is an 
isomorphism
\begin{equation} \label{eqn1}
\rho : M \iso \opn{RHom}_{B \otimes_{\K} B}(B, M \otimes_{\K} M)
\end{equation}
in $\msf{D}(\cat{Mod} B)$.

Now suppose $A$ is any ring. 
Trying to write $A$ instead of $\K$ in
formula (\ref{eqn1}) does not make sense:
instead of $M \otimes_{A} M$ we must take the derived tensor 
product
$M \otimes^{\mrm{L}}_{A} M$; but then there is no obvious way to 
make $M \otimes^{\mrm{L}}_{A} M$ into a complex of 
$B \otimes_{A} B$ -modules. 

The problem is torsion: $B$ might fail to be a flat $A$-algebra. 
This is where {\em differential graded algebras} (DG algebras) 
enter the picture.

A DG algebra is a graded ring
$\til{A} = \boplus_{i \in \mbb{Z}} \til{A}^i$, 
together with a graded derivation
$\mrm{d} : \til{A} \to \til{A}$ of degree $1$, satisfying 
$\mrm{d} \circ \mrm{d} = 0$. 

A DG algebra quasi-isomorphism is a homomorphism 
$f : \til{A} \to \til{B}$ respecting degrees, multiplications and 
differentials, and such that 
$\mrm{H}(f) : \mrm{H} \til{A} \to \mrm{H} \til{B}$
is an isomorphism (of graded algebras).

We shall only consider {\em super-commutative non-positive} DG 
algebras. Super-commutative means that 
$a b = (-1)^{i j} b a$ and $c^2 = 0$
for all $a \in \til{A}^i$, $b \in \til{A}^j$ and 
$c \in \til{A}^{2i + 1}$. Non-positive means that
$\til{A} = \boplus_{i \leq 0} \til{A}^i$.

We view a ring $A$ as a DG algebra concentrated in degree $0$. 
Given a DG algebra homomorphism $A \to \til{A}$ we say that 
$\til{A}$ is a DG $A$-algebra.

Let $A$ be a ring. 
A {\em semi-free} DG $A$-algebra is a DG $A$-algebra 
$\til{A}$, such that after forgetting the differential $\til{A}$
is isomorphic, as graded $A$-algebra, to a super-polynomial algebra
on some graded set of variables. 

\begin{dfn}
Let $A$ be a ring and $B$ an $A$-algebra. 
A {\em semi-free DG algebra resolution of $B$ relative to $A$} 
is a quasi-isomorphism  $\til{B} \to B$ of DG 
$A$-algebras, where $\til{B}$ is a 
semi-free DG $A$-algebra.
\end{dfn}

Such resolutions always exist, and they are unique up to 
quasi-isomorphism.  

\begin{exa}
Take $A = \mbb{Z}$ and $B = \mbb{Z} / (6)$.
Define $\til{B}$ to be the super-polynomial algebra
$A[\xi]$ on the variable $\xi$ of degree $-1$. 
So $\til{B} = A \oplus A \xi$
as free graded $A$-module, and $\xi^2 = 0$.
Let $\mrm{d}(\xi) := 6$. Then $\til{B} \to B$ is a 
semi-free DG algebra resolution of $B$ relative to $A$.
\end{exa}

For a DG algebra $\til{A}$ one has the category $\cat{DGMod} \til{A}$ 
of DG $\til{A}$-modules. It is analogous to the category of 
complexes of 
modules over a ring, and by a similar process of inverting 
quasi-isomorphisms we obtain the derived category 
$\til{\cat{D}}(\cat{DGMod} \til{A})$; see \cite{Ke},
\cite{Hi}. 

For a ring $A$ (i.e.\ a DG algebra concentrated in degree 
$0$) we have
\[ \til{\cat{D}}(\cat{DGMod} A) = \msf{D}(\cat{Mod} A) , \]
the usual derived category.

It is possible to derive functors of DG modules, again in analogy 
to $\msf{D}(\cat{Mod} A)$.
An added feature is that for a 
quasi-isomorphism $\til{A} \to \til{B}$, the restriction of scalars functor
\[ \til{\cat{D}}(\cat{DGMod} \til{B}) \to
\til{\cat{D}}(\cat{DGMod} \til{A}) \]
is an equivalence.

Getting back to our original problem, suppose $A$ is a ring and 
$B$ is an $A$-algebra. Choose a semi-free DG algebra resolution 
$\til{B} \to B$ relative to $A$. 
For $M \in \msf{D}(\cat{Mod} B)$ define
\[ \opn{Sq}_{B / A} M :=
\mrm{RHom}_{\til{B} \otimes_{A} \til{B}}(B, 
M \otimes^{\mrm{L}}_{A} M) \]
in $\msf{D}(\cat{Mod} B)$.

\begin{thm}\tup{(\cite{YZ4})}
The functor
\[ \opn{Sq}_{B / A} : \msf{D}(\cat{Mod} B) \to
\msf{D}(\cat{Mod} B) \]
is independent of the resolution $\til{B} \to B$.
\end{thm}

The functor $\opn{Sq}_{B / A}$, called the {\em squaring 
operation}, is nonlinear. In fact, given
a morphism $\phi : M \to M$ in $\msf{D}(\cat{Mod} B)$ and an 
element $b \in B$ one has
\begin{equation} \label{eqn2}
\opn{Sq}_{B / A}(b \phi) = b^2 \opn{Sq}_{B / A}(\phi) 
\end{equation}
in
\[ \opn{Hom}_{\msf{D}(\cat{Mod} B)}(\opn{Sq}_{B / A} M,
\opn{Sq}_{B / A} M) . \]

\begin{dfn}
Let $B$ be a noetherian $A$-algebra, and let
$M$ be a complex in 
$\msf{D}^{\mrm{b}}_{\mrm{f}}(\cat{Mod} B)$ that has finite flat 
dimension over $A$. Assume
\[ \rho : M \iso \opn{Sq}_{B / A} M \]
is an isomorphism in $\msf{D}(\cat{Mod} B)$. Then the pair
$(M, \rho)$ is called a {\em rigid complex over $B$ relative to 
$A$}.
\end{dfn}

\begin{dfn} 
Say $(M, \rho)$ and $(N, \sigma)$ are rigid complexes over $B$ 
relative to $A$. A morphism $\phi : M \to N$ in
$\msf{D}(\cat{Mod} B)$ is called a {\em 
rigid morphism relative to $A$} if the diagram
\[ \begin{CD}
M @>{\rho}>> \opn{Sq}_{B / A} M \\
@V{\phi}VV @VV{\opn{Sq}_{B / A} (\phi)}V \\
N  @>{\sigma}>> \opn{Sq}_{B / A} N
\end{CD} \]\
is commutative. 
\end{dfn}

We denote by 
$\msf{D}^{\mrm{b}}_{\mrm{f}}(\cat{Mod} B)_{\mrm{rig} / A}$
the category of rigid complexes over $B$ relative to $A$. 

\begin{exa} \label{exa1}
Take $M = B = A$. Then
\[ \opn{Sq}_{A / A} A = 
\mrm{RHom}_{A \otimes_{A} A}(A, A \otimes_{A} A) = A , \]
and we interpret this as a rigidifying isomorphism
\[ \rho^{\mrm{tau}} : A \iso \opn{Sq}_{A / A} A . \]
The {\em tautological rigid complex} is
\[ (A, \rho^{\mrm{tau}}) \in 
\msf{D}^{\mrm{b}}_{\mrm{f}}(\cat{Mod} A)_{\mrm{rig} / A} . \]
\end{exa}

\section{Properties of Rigid Complexes}

The first property of rigid complexes explains their name.

\begin{thm} \tup{(\cite{YZ4})}
Let $A$ be a ring, $B$ a noetherian $A$-algebra, and 
\[ (M, \rho) \in 
\msf{D}^{\mrm{b}}_{\mrm{f}}(\cat{Mod} B)_{\mrm{rig} / A} . \]
Assume the canonical ring homomorphism 
\[ B \to \opn{Hom}_{\msf{D}(\cat{Mod} B)}(M, M) \]
is bijective. Then the only automorphism of $(M, \rho)$
in 
$\msf{D}^{\mrm{b}}_{\mrm{f}}(\cat{Mod} B)_{\mrm{rig} / A}$
is the identity $\bsym{1}_M$. 
\end{thm}

The proof is very easy: an automorphism $\phi$ of $M$ 
has to be of the form $\phi = b\, \bsym{1}_M$ for some invertible 
element $b \in B$. If $\phi$ is rigid then $b = b^2$
(cf.\ formula (\ref{eqn2})), and hence 
$b = 1$. 

We find it convenient to denote ring homomorphisms by $f^*$ etc. 
Thus a ring homomorphism $f^* : A \to B$ corresponds to the 
morphism of schemes
\[ f : \opn{Spec} B \to \opn{Spec} A . \]

Let $A$ be a noetherian ring. 
Recall that an $A$-algebra $B$ is called {\em essentially finite type} 
if it is a localization of some finitely generated $A$-algebra. 
We say that $B$ is {\em essentially smooth}
(resp.\ {\em essentially \'etale}) over $A$ if it is essentially 
finite type and formally smooth (resp.\ formally \'etale).

\begin{exa}
If $A'$ is a localization of $A$ then $A \to A'$ is essentially 
\'etale. If $B = A[t_1, \ldots, t_n]$ is a polynomial algebra then 
$A \to B$ is smooth, and hence also essentially smooth.
\end{exa}

Let $A$ be a noetherian ring and $f ^* : A \to B$ an 
essentially smooth homomorphism. Then $\Omega^1_{B / A}$ is a 
finitely generated projective $B$-module. Let
\[ \opn{Spec} B = \coprod_i \opn{Spec} B_i \]
be the decomposition into connected components, and for every $i$ 
let $n_i$ be the rank of $\Omega^1_{B_i / A}$. We define a functor
\[ f^{\sharp} : \msf{D}(\cat{Mod} A) \to \msf{D}(\cat{Mod} B) \]
by
\[ f^{\sharp} M := \bigoplus_i \, 
\Omega^{n_i}_{B_i / A}[n_i] \otimes_A M . \]

Recall that a ring homomorphism $f^* : A \to B$ is called {\em finite} 
if $B$ is a finitely generated $A$-module. Given such a finite 
homomorphism we define a functor
\[ f^{\flat} : \msf{D}(\cat{Mod} A) \to \msf{D}(\cat{Mod} B) \]
by
\[ f^{\flat} M := \opn{RHom}_A(B, M) . \]

\begin{thm}\tup{(\cite{YZ4})} \label{thm0.1}
Let $A$ be a noetherian ring, let $B, C$ be essentially finite 
type $A$-algebras, let $f^* : B \to C$ be an $A$-algebra
homomorphism, and let 
\[ (M, \rho) \in 
\msf{D}^{\mrm{b}}_{\mrm{f}}(\cat{Mod} B)_{\mrm{rig} / A} . \]
\begin{enumerate}
\rmitem{1} If $f^*$ is finite and $f^{\flat} M$ has finite flat 
dimension over $A$, then $f^{\flat} M$ has an induced rigidifying isomorphism 
\[ f^{\flat}(\rho) : f^{\flat} M  \iso \opn{Sq}_{C / A}
f^{\flat} M . \]
The assignment 
\[ (M, \rho) \mapsto f^{\flat} (M, \rho) :=
\bigl( f^{\flat}(\rho), f^{\flat} M \bigr) \]
is functorial. 
\rmitem{2} If $f^*$ is essentially smooth
then $f^{\sharp} M$ has an induced rigidifying isomorphism 
\[ f^{\sharp}(\rho) : f^{\sharp} M  \iso \opn{Sq}_{C / A}
f^{\sharp} M . \]
The assignment 
\[ (M, \rho) \mapsto f^{\sharp} (M, \rho) :=
\bigl( f^{\sharp}(\rho), f^{\sharp} M \bigr) \]
is functorial. 
\end{enumerate}
\end{thm}

\section{Rigid Dualizing Complexes}

Let $\K$ be a regular noetherian ring of finite Krull dimension. 
We denote by $\cat{EFTAlg} / \K$ the category of essentially finite type 
$\K$-algebras. 

\begin{dfn}
A {\em rigid dualizing complex} over $A$ relative to $\K$ is a 
rigid complex $(R_A, \rho_A)$, such that $R_A$ is a dualizing 
complex. 
\end{dfn}

\begin{thm} \tup{(\cite{YZ5})} \label{thm4}
Let $\K$ be a regular finite dimensional noetherian ring, and
let $A$ be an essentially finite type $\K$-algebra.
\begin{enumerate}
\rmitem{1} The algebra $A$ has a rigid dualizing complex 
$(R_A, \rho_A)$, which is unique up to a unique rigid isomorphism.
\rmitem{2} Given a finite homomorphism $f^* : A \to B$, there is a 
unique rigid isomorphism
$f^{\flat}(R_A, \rho_A) \iso (R_B, \rho_B)$.
\rmitem{3} Given an essentially smooth homomorphism $f^* : A \to B$ , 
there is a unique rigid isomorphism
$f^{\sharp}(R_A, \rho_A) \iso (R_B, \rho_B)$.
\end{enumerate}
\end{thm}

Here is how the rigid dualizing complex $(R_A, \rho_A)$ is 
obtained. We begin with the tautological 
rigid complex 
\[ (\K, \rho^{\mrm{tau}}) \in
\msf{D}^{\mrm{b}}_{\mrm{f}}(\cat{Mod} \K)_{\mrm{rig} / \K} , \]
which is dualizing (cf.\ Examples \ref{exa0} and \ref{exa1}). 
Now the structural homomorphism $\K \to A$ can be factored into
\[ \K \xar{f^*} B \xar{g^*} C \xar{h^*} A , \]
where $f^*$ is essentially smooth ($B$ is a polynomial algebra over $\K$);
$g^*$ is finite (a surjection); 
and $h^*$ is also essentially smooth (a localization). 

It is not hard to check (see \cite[Chapter V]{RD}) that each of the complexes
$f^{\sharp} \K$, $g^{\flat} f^{\sharp} \K$ and
$h^{\sharp}\, g^{\flat} f^{\sharp} \K$ 
is dualizing over the respective ring. In particular, $g^{\flat} f^{\sharp} \K$
has bounded cohomology, and hence it has finite flat dimension over 
$\K$. According to Theorem \ref{thm0.1} we then have a rigid complex
\[ (R_A, \rho_A) := h^{\sharp}\, g^{\flat}\, f^{\sharp}
(\K, \rho^{\mrm{tau}}) \in
\msf{D}^{\mrm{b}}_{\mrm{f}}(\cat{Mod} A)_{\mrm{rig} / \K} . \]

\begin{dfn}
Given a homomorphism $f^* : A \to B$ in $\cat{EFTAlg} / \K$,
define the {\em twisted inverse image functor}
\[ f^! : \msf{D}^{+}_{\mrm{f}}(\cat{Mod} A) \to
\msf{D}^{+}_{\mrm{f}}(\cat{Mod} B) \]
by the formula
\[ f^! M :=  
\opn{RHom}_B \bigl( B \otimes^{\mrm{L}}_{A} 
\opn{RHom}_A(M, R_A), R_B \bigr) . \]
\end{dfn}

It is easy to show that the assignment $f^* \mapsto f^!$ 
is a $2$-functor from the category $\cat{EFTAlg} / \K$ to the $2$-category
$\cat{Cat}$ of all categories. 
Moreover, using Theorem \ref{thm4} one can show that this operation has very 
good properties. For instance, when $f^*$ is finite, 
then there is a functorial nondegenerate trace morphism
\[ \opn{Tr}_f :  f^{!} M \to M . \]

\begin{rem} \label{rem1} 
According to Grothendieck's duality theory in \cite{RD}, if 
$f : X \to Y$ is a finite type morphism between noetherian schemes,
and if $Y$ has a dualizing complex, then there is a functor
\[ f^{! (\mrm{G})} : \msf{D}^{+}_{\mrm{c}}(\cat{Mod} \mcal{O}_Y) \to
\msf{D}^{+}_{\mrm{c}}(\cat{Mod} \mcal{O}_X) , \]
with many good properties. 

Let $\cat{FTAlg} / \K$ be the category of finite type $\K$-algebras. 
By restricting attention to affine schemes, the results of \cite{RD}
give rise to a $2$-functor $f^* \mapsto f^{! (\mrm{G})}$
from $\cat{FTAlg} / \K$ to $\cat{Cat}$. 
It is not hard to show that the $2$-functor 
$f^* \mapsto f^{! (\mrm{G})}$ is isomorphic to our 
$2$-functor $f^* \mapsto f^{!}$; see \cite[Theorem 4.10]{YZ5}.

It should be noted that our construction works in the slightly bigger
category $\cat{EFTAlg} / \K$. It also has the advantage of
being local; whereas in \cite{RD} some of the results require that
morphisms between affine schemes be compactified.
\end{rem}

\section{Rigid Complexes and CM Homomorphisms}

In this final section we discuss the relation between rigid 
complexes and Cohen-Macaulay homomorphisms. 

\begin{dfn}
A ring $A$ is called {\em tractable} if there is an 
essentially finite type homomorphism $\K \to A$, for some 
regular noetherian ring of finite Krull dimension $\K$.
\end{dfn}

Such a homomorphism $\K \to A$ is called a {\em traction} for $A$. 
It is not part of the structure -- the ring $A$ does come with
any preferred traction. 
``Most commutative noetherian rings we know'' are tractable.

Given a traction $\K \to A$ we denote by $R_{A / \K}$ the rigid
dualizing complex of $A$ relative to $\K$; cf.\ Theorem \ref{thm4}. 
(The rigidifying isomorphism $\rho_{A / \K}$ is implicit.)

Recall that a noetherian ring $A$ is called {\em Cohen-Macaulay} 
(resp.\ {\em Gorenstein}) if all its local rings $A_{\mfrak{p}}$,
$\mfrak{p} \in \opn{Spec} A$, are Cohen-Macaulay (resp.\ Gorenstein) 
local rings. The implications are regular $\Rightarrow$ Gorenstein
$\Rightarrow$ Cohen-Macaulay.

Let $f^* : A \to B$ be a ring homomorphism. For 
$\mfrak{p} \in \opn{Spec} A$ let 
$\bsym{k}(\mfrak{p}) := (A / \mfrak{p})_{\mfrak{p}}$, the 
residue field. The fiber of $f^*$ above $\mfrak{p}$ is the 
$\bsym{k}(\mfrak{p})$-algebra
$B \otimes_A \bsym{k}(\mfrak{p})$. 
Now assume $f^*$ is an essentially finite type flat homomorphism. 
If all the fibers of $f^*$ are Cohen-Macaulay (resp.\ Gorenstein) 
rings, then we call $f^*$ an {\em essentially Cohen-Macaulay} (resp.\ 
{\em essentially Gorenstein}) homomorphism.

\begin{thm}\tup{(\cite{Ye2})} \label{thm6}
Let $A$ be a tractable ring, and let $f^* : A \to B$ be homomorphism
which is of essentially finite type and of finite flat dimension.
Then there exists a rigid complex $R_{B/A}$ over $B$ relative to $A$,
unique up to a unique rigid isomorphism, with the following property:
\begin{enumerate}
\rmitem{*} Let $\K \to A$ be some traction. Then 
\[ R_{A / \K} \otimes^{\mrm{L}}_{A} R_{B / A} \cong
R_{B / \K} \]
in $\msf{D}(\cat{Mod} B)$. 
\end{enumerate}
\end{thm}

Condition (*) implies that the support of the complex $R_{B / A}$ is 
$\opn{Spec} B$.
One can prove that 
\[ f^! M \cong R_{B / A} \otimes^{\mrm{L}}_{A} M \]
for $M \in \msf{D}^{\mrm{b}}_{\mrm{f}}(\cat{Mod} A)$.

If the ring $A$ is Gorenstein, then $R_{A / \K}$ is a shift of an
invertible $A$-module. Hence:

\begin{cor} \label{cor2}
Assume that in Theorem \tup{\ref{thm6}} the ring $A$ is Gorenstein.
Then 
$R_{B / A}$ is a dualizing complex over $B$ 
\end{cor}

The rigid complex $R_{B / A}$ allows us to characterize
Cohen-Macaulay 
homomorphisms, as follows.

\begin{thm}\tup{(\cite{Ye2})} \label{thm5}
Let $A$ be a tractable ring, and let $f^* : A \to B$ be an essentially finite 
type flat homomorphism. Then the following conditions are equivalent:
\begin{enumerate}
\rmitem{i} $f^*$ is an essentially Cohen-Macaulay homomorphism.
\rmitem{ii} Let 
\[ \opn{Spec} B = \coprod_i \opn{Spec} B_i \]
be the decomposition into connected components. Then for any $i$ 
there is a finitely generated $B_i$-module
$\bsym{\omega}_{B_i / A}$, which is flat over $A$,
and an integer $n_i$, such that 
\[ R_{B / A} \cong \bigoplus_i \bsym{\omega}_{B_i / A}[n_i] \]
in $\msf{D}(\cat{Mod} B)$. 
\end{enumerate}
\end{thm}

The module 
\[ \bsym{\omega}_{B / A} := \bigoplus_i \bsym{\omega}_{B_i / A} \]
is called the {\em relative dualizing module} of
$f^* : A \to B$. Note that the complex 
$\boplus_i \bsym{\omega}_{B_i / A}[n_i]$ is rigid, but in general 
it is not a dualizing complex over $B$. Still the fibers of 
$\boplus_i \bsym{\omega}_{B_i / A}[n_i]$ are dualizing 
complexes -- this can be seen by taking $A' = \bsym{k}(\mfrak{p})$ in the 
next result, and using Corollary \ref{cor2}

Here is a ``rigid'' version of Conrad's base change theorem 
\cite{Co}.

\begin{thm}\tup{(\cite{Ye2})} 
Let 
\[ \begin{CD}
A @>{}>> B \\
@V{}VV @VV{}V \\
A'  @>{}>> B'
\end{CD} \]
be a cartesian diagram of rings, i.e.\ 
\[ B' \cong A' \otimes_A B , \]
with $A$ and $A'$ tractable rings. Assume 
$A \to B$ is an essentially Cohen-Macaulay homomorphism. \tup{(}There isn't any  
restriction on the homomorphism $A \to A'$.\tup{)} Then:
\begin{enumerate}
\rmitem{1} $A' \to B'$ is an essentially Cohen-Macaulay homomorphism.
\rmitem{2} There is a unique isomorphism of $B'$-modules
\[ \bsym{\omega}_{B' / A'} \cong A' \otimes_A \bsym{\omega}_{B / A}  \]
which respects rigidity.
\end{enumerate}
\end{thm}

 From this we can easily deduce the next result.

\begin{cor} \label{cor1}
Let $A$ be a tractable ring, and let $f^* : A \to B$ be an essentially 
Cohen-Macaulay homomorphism. Then the following conditions are equivalent:
\begin{enumerate}
\rmitem{i} $f^*$ is an essentially Gorenstein homomorphism.
\rmitem{ii} $\bsym{\omega}_{B / A}$ is an invertible $B$-module.
\end{enumerate}
\end{cor}

\begin{rem}
The recent paper \cite{AI} contains results similar to Theorem 
\ref{thm5} and Corollary \ref{cor1}, obtained by 
different methods, and without the requirement that $A$ is tractable.
\end{rem}

\end{document}